\documentclass[11pt]{amsart}

\usepackage{amsfonts}
\usepackage{amssymb}
\usepackage{amsmath}
\usepackage{amsthm}
\usepackage{amsopn}
\usepackage{graphicx}
\usepackage{times}
\usepackage[normalem]{ulem}

\newcommand{\F}{\mathbb F}
\newcommand{\N}{\mathbb N}
\newcommand{\Z}{\mathbb Z}
\newcommand{\R}{\mathbb R}
\newcommand{\gldz}{\mbox{\rm GL}_d(\Z)}

\DeclareMathOperator{\vol}{vol}
\DeclareMathOperator{\conv}{conv}
\DeclareMathOperator{\mv}{mv}
\DeclareMathOperator{\dt}{dt}

\newtheorem{theorem}{Theorem}

\newtheorem*{remark}{Remark}

\newtheorem*{proposition}{Proposition}
\newtheorem{corollary}{Corollary}

\renewcommand{\vec}[1]{\boldsymbol{#1}}

\input cyracc.def
\font\tencyr=wncyr10
\def\cyr{\tencyr\cyracc}

\begin{document}

\title{Lattice Delone simplices with super-exponential volume}

\author{Francisco Santos}
\address{Dept. Matem\'aticas, Estad\'{\i}stica y
         Computaci\'on, Univ. de Cantabria, Spain}
\email{santosf@unican.es}
\thanks{The first author was partially 
        supported by the Spanish Ministry of Science and Education
        under grant MTM2005-08618-C02-02}

\author{Achill Sch\"urmann}
\thanks{The second and third author were partially 
        supported by the ``Deutsche Forschungsgemeinschaft'' (DFG) 
        under grant SCHU 1503/4-1.}

\address{Mathematics Department, University of Magdeburg, 
         39106 Magdeburg, Germany}
\email{achill@math.uni-magdeburg.de}

\author{Frank Vallentin}
\thanks{The third author was partially supported by the Edmund Landau Center for 
                 Research in Mathematical Analysis and Related Areas, 
                 sponsored by the Minerva Foundation (Germany).}
\address{Einstein Institute of Mathematics, The Hebrew
         University of Jerusalem, Jerusalem, 91904, Israel}
\email{frank.vallentin@gmail.com}

\begin{abstract}
In this short note we give a construction of an infinite series of 
Delone simplices whose relative volume grows super-exponentially with 
their dimension. This dramatically improves the previous best lower bound,
which was linear.
\end{abstract}

\maketitle

\section*{Background}

Consider the Euclidean space $\R^d$ with norm $\|\cdot\|$ and a
discrete subset $\Lambda \subset \R^d$.  A $d$-dimensional polytope $L
= \conv\{\vec{v}_0, \ldots, \vec{v}_n\}$ with $\vec{v}_0, \ldots,
\vec{v}_n\in \Lambda$ is called a \textit{Delone polytope} of
$\Lambda$, if there exists an empty sphere $S$ with 
$S\cap \Lambda =\{\vec{v}_0, \ldots, \vec{v}_n\}$. That is, 
if there is 
a center $\vec{c}\in\R^d$ and a radius
$r>0$ such that $\|\vec{v}_i-\vec{c}\| = r$ for $i=0,\dots,n$, and
$\|\vec{v}-\vec{c}\|>r$ for the remaining
$\vec{v}\in\Lambda\setminus\{\vec{v}_0, \ldots, \vec{v}_n\}$.  If the
Delone polytope is a simplex, hence $n=d$, we speak of a
\textit{Delone simplex}.  Over the past decades Delone polytopes
experienced a renaissance in applications like computer graphics and computational geometry, 
where they are traditionally called Delaunay polytopes due to the French
transcription of {\cyr Delone}.

Historically, Delone polytopes received attention in the study of
positive definite quadratic forms (PQFs), in particular in a reduction
theory due to Voronoi (cf. \cite{voronoi-1908}).  To every $d$-ary PQF
$Q$, a point lattice $\Lambda=A\Z^d$ with $Q=A^tA$ is associated;
$\Lambda$ is a discrete set and uniquely determined by $Q$ up to
orthogonal transformations.  A $d$-dimensional polytope $L =
\conv\{\vec{v}_0, \ldots, \vec{v}_n\}$, with $\vec{v}_0, \ldots,
\vec{v}_n\in \Z^d$, is a Delone polytope of the lattice
$\Lambda$ (and also called a Delone polytope of ~$Q$), if and only if
there exists a $\vec{c} \in \R^d$ and $r > 0$ with $Q[\vec{v} -
\vec{c}] = \|A(\vec{v} - \vec{c})\|^2 \geq r^2$ for all
$\vec{v}\in\Z^d$ and with equality if and only if $\vec{v}=\vec{v}_i$
for $i\in\{0,\dots,n\}$.  The set of all Delone polytopes forms a
periodic face-to-face tiling of $\R^d$. It is called the
\textit{Delone subdivision} of~$Q$.  If all Delone polytopes are
simplices we speak of a \textit{Delone triangulation}. Delone
subdivisions form a poset with respect to refinement, in which
triangulations are maximally refined elements.  

Two Delone polytopes
$L$ and $L'$ are \textit{unimodularly equivalent}, if $L=UL'+\vec{t}$
for some unimodular transformation $U\in\gldz$ and a translation
vector $\vec{t}\in\Z^d$. Voronoi \cite{voronoi-1908} showed that, up to
unimodular equivalence, there exist only finitely many Delone
subdivisions in each dimension $d$.  
Delone simplices may be considered as ``building
blocks'' in this theory and therefore their classification is of
particular interest.  For more information on the classical theory 
we refer the interested reader to \cite{sv-2004} or the original 
works of Voronoi \cite{voronoi-1908} and Delone \cite{delone-1937}.

\section*{Bounds for relative volume}

Let $\Lambda\subset \R^d$ be a lattice.  A lattice simplex
$\conv\{\vec{v}_0,\dots,\vec{v}_n\}$ is called unimodular if
$\{\vec{v}_0,\dots,\vec{v}_n\}$ is an affine basis of $\Lambda$. All
unimodular simplices are unimodularly equivalent and, in particular,
have the same volume.  The \textit{relative volume} (or
\textit{normalized volume}) of a lattice simplex $L$ is the volume of
$L$ divided by that of a unimodular simplex of $\Lambda$,
so that the relative volume equals $\vol(L) \cdot d!$.
Equivalently, it equals the index in $\Lambda$ of
the sublattice affinely spanned by the vertices of $L$.

Clearly, the relative volume is 
an invariant with respect to unimodular transformations.
Hence, in order to classify
possible Delone simplices up to unimodular equivalence, a first
question, already raised by Delone in \cite{delone-1937}, is what is the
maximum relative volume $\mv(d)$ of $d$-dimensional Delone simplices. The sequence $\mv(d)$ is (weakly) increasing, since
from a Delone simplex  for a lattice $\Lambda\subset \R^d$
one can easily construct another of the same relative volume for the 
lattice $\Lambda\times \Z\subset \R^{d+1}$.

Voronoi knew that up to dimension $d = 4$ all Delone simplices have
relative volume $1$, while there are Delone simplices of volume 2 in
$d = 5$.  In \cite{baranovskii-1973} Baranovskii proved $\mv(5) = 2$,
and later Baranovskii and Ryshkov \cite{br-1998} proved $\mv(6) =
3$. Dutour classified all $6$-dimensional Delone polytopes in
\cite{dutour-2004}. 

Ryshkov \cite{ryshkov-1973c} 
was the first who proved that relative volumes of Delone 
simplices are not bounded when the dimension goes to infinity. 
More precisely, for every $k\in\N$ he constructed Delone simplices of relative volume $k$ in dimension $2k+1$,
establishing that $\mv(d) \geq
\left\lfloor\frac{d-1}{2}\right\rfloor$. This was recently
improved to the still linear lower bound 
$\mv(d) \geq d - 3$ by Erdahl and
Rybnikov \cite{er-2002}. 

In this note we prove the following two lower bounds on $\mv(d)$:

\begin{theorem}
For every pair of dimensions $d_1$ and $d_2$,
\[
\mv(d_1+d_2) \ge \mv(d_1)\mv(d_2).
\]
\end{theorem}

\begin{theorem}
For every dimension of the form $d=2^n-1$,
$$\displaystyle \mv(d) \geq (d+1)^{(d+3)/2} / 4^d.$$
\end{theorem}

Theorem 1 immediately implies exponential lower bounds 
on $\mv(d)$. For example, $\mv(5)=2$ gives
$
\mv(d)\ge 2^{\lfloor d/5 \rfloor}\sim 1.1487^d.
$
Even better, it is known that $\mv(24)\ge 20480$. Indeed, the Delone
subdivision of the Leech lattice, which was
determined by Borcherds, Conway, Queen, Parker and Sloane
\cite[Chapter~25]{cs-1988}, contains simplices of
relative volume $20480$ (the simplex denoted $a^{24}_1a_1$
in their classification). Therefore we obtain:
\begin{corollary}
\[
\mv(d)\geq 20480^{\left\lfloor d/24\right\rfloor}  
\sim 1.5123^d.
\]
\end{corollary}

Theorem 2 gives a much better lower bound asymptotically:

\begin{corollary}
\[
\log(\mv(d))\in \Theta(d\log d).
\]
\end{corollary}

\begin{proof}
For the lower bound, let $2^n$ be the largest power of two that is
smaller or equal to $d+1$ (so that $2^n\ge (d+1)/2$). Theorem 2,
together with the monotonicity of $\mv(d)$, gives:
\[
\mv(d)\ge \mv(2^n-1)\ge 
\frac{(2^n)^{(2^n+2)/2}}{4^{2^n-1}}
\in 2^{\Theta\left(n2^n \right)} =
 2^{\Theta\left(d \log d\right)}.
\]

For the upper bound we use the following argument, which is
Lemma 14.2.5 in~\cite{dl-1997} (attributed to
L.~Lovasz): Given a Delone simplex $L$
of some PQF, the volume of the centrally symmetric difference body
$L-L=\{\vec{v}-\vec{v'} : \vec{v},\vec{v}'\in L \}$ is
$\vol(L-L)=\binom{2d}{d}\vol(L)$ (see \cite{rs-1957}).  The polytope
$L-L$ does not contain elements of $L\setminus\{0\}$ in its
interior.
Thus, by Minkowski's fundamental theorem (see
\cite{gl-1987}, \S 5 Theorem~1) we know that $\vol(L-L)\leq
2^d$. Putting things together we get 
\[
\mv(d) \le  \frac{2^d d!}{\binom{2d}{d}}
\sim \sqrt{2}\pi d\left(\frac{d}{2e}\right)^d.
\]
\end{proof}

More precisely, the arguments in this proof say that
\[
\frac{1}{4}\le \liminf {\frac{\log \mv(d)}{d\log d}
\le 
\limsup \frac{\log \mv(d)}{d\log d}
 \le 1}.
\]
Similarly, Theorem 2 directly implies 
\[
\frac{1}{2}\le 
\limsup \frac{\log \mv(d)}{d\log d}.
\]
We do not know whether 
$ \lim \frac{\log \mv(d)}{d\log d}$ exists.


%
%
%

\section*{Proof of Theorem 1}

Theorem 1 follows from the fact that an orthogonal product
of simplices decomposes into simplices with relative volume being the
product of the individual relative volumes. Let us be more precise:
Let $L_1 = \conv\{\vec{v}_0, \ldots, \vec{v}_{d_1}\} \subseteq \R^{d_1}$ 
be a Delone simplex of the 
lattice $\Lambda_1$ with relative volume $\mv(d_1)$, and let $L_2 = \conv\{\vec{w}_0,
\ldots, \vec{w}_{d_2}\} \subseteq \R^{d_2}$ be a Delone simplex of the
lattice $\Lambda_2$ 
with relative volume $\mv(d_2)$.  Then, the direct product $L_1
\times L_2 = \conv\{(\vec{v}_i, \vec{w}_j) : i = 0, \ldots, d_1,\; j =
0, \ldots, d_2\} \subseteq \R^{d_1} \times \R^{d_2}$ is a $(d_1 +
d_2)$-dimensional Delone polytope of the lattice $\Lambda= \Lambda_1\times
\Lambda_2$. Let $\Lambda'_1$, $\Lambda'_2$ and $\Lambda'=\Lambda'_1\times \Lambda'_2$ denote the lattices affinely generated by the vertices of $L_1$,
$L_2$ and $L_1\times L_2$, respectively. 
By the classical theory of Voronoi (see \cite{sv-2004}, Proposition
5.1 and Proposition 5.4) we know that by a suitable infinitesimal
change of the PQF $Q$ that induces $\Lambda$ the Delone polytope $L_1 \times L_2$ is triangulated into Delone simplices, which hence have vertices in the sublattice $\Lambda'$
(more precisely, in a perturbation of it). Since the index of $\Lambda'$ in $\Lambda$
is precisely $\mv(d_1) \mv(d_2)$, these Delone simplices have
relative volume at least that number.

\begin{remark}\rm
The product of two simplices $L_1$ and $L_2$ is a
\textit{totally unimodular} polytope, meaning that all the simplices
spanned by a subset of its vertices have the same volume. In
particular, the Delone simplices that we obtain in the last step of
the proof have relative volume \textit{exactly} $\mv(d_1) \mv(d_2)$
(cf., for example, \cite{haiman-1991}). As a consequence, every triangulation
of $L_1 \times L_2$ consists of exactly $\binom{d_1 + d_2}{d_1}$ simplices.
\end{remark}

\section*{Proof of Theorem 2}

Remember that a Hadamard matrix of order $d$ is a $d \times d$ matrix 
with elements $+1$ and $-1$ in which distinct columns are orthogonal. 
Hadamard matrices exist at least whenever $d$ is a power of two and conjecturally
whenever $d$ is a multiple of four (cf. \cite{cs-1988}, Chapter~3,
2.13; for a recent survey on Hadamard matrices see~\cite{ek05}).

By
multiplying columns with $\pm 1$ we normalize columns of a Hadamard
matrix $H$ so that the
entries of the first row are all $+1$. Then, from $H$ we get the $(d-1) \times d$
matrix $\tilde H$ with elements $0$ and $1$ by deleting the
first row and replacing $+1$ by $0$, and $-1$ by $1$. The columns
of $\tilde H$ form the vertex set of a regular $(d-1)$ simplex with
edge length $\sqrt{d/2}$ and volume $\det(H)/2^{d-1} = d^{d/2}/2^{d-1}$.
This is the maximum
volume of a simplex contained in the $(d-1)$ unit cube
$[0,1]^{d-1}$. 
Such a simplex is called a $(d-1)$ Hadamard simplex.

A standard construction for Hadamard matrices of order $2^n$ is as follows
\[
H_1 = (1),\quad
H_{2^{n}} =
\begin{pmatrix}
H_{2^{n-1}} & H_{2^{n-1}}\\
-H_{2^{n-1}} & H_{2^{n-1}}
\end{pmatrix}.
\]
Let $d=2^n-1$ and consider the matrix $\tilde H_{2^n}$ associated
to the particular Hadamard matrix $H_{2^n}$ constructed this way.
It is well known, and easy to see, that the columns of $\tilde H_{2^n}$ 
form an $n$-dimensional subspace
of ${\F_2}^{d}$. This subspace is also known
as the $[d,n]$ binary simplex code or the dual of the $[d,
d-n]$ binary Hamming code. Linearity implies that the following
is a $d$-dimensional sublattice of $\Z^{d}$:
\[
\Lambda(\tilde H_{2^n}) = \{(v_1, \ldots, v_{d}) \in \Z^{d} : (v_1
\operatorname{mod}
2, \ldots, v_{d} \operatorname{mod} 2) \in \tilde H_{2^n}\}.
\] 
(This procedure is ``Construction A'' in \cite{cs-1988}, Chapters~5 and~7).
The lattice $\Lambda(\tilde H_{2^n})$ has determinant
$2^{d-n}$ in $\Z^{d}$, since it contains $2^n$ of the 
$2^{d}$ vertices in every lattice unit cube.
Moreover, the Hadamard simplex defined by  $\tilde H_{2^n}$
is a Delone simplex in $\Lambda(\tilde H_{2^n})$
because the sphere around the unit cube does
not contain more lattice points than the columns of $\tilde H_{2^n}$.
 The relative volume of the Hadamard simplex is
\[
\frac{(d+1)^{(d+1)/2}/2^{d}}{2^{d-n}} =\frac{(d+1)^{(d+3)/2}}{4^{d}}.
\]

\section*{Remarks on related problems}

\begin{enumerate}

\item 
The determination of $\mv(d)$ is clearly related (but not equivalent) to 
that of what is the minimal number
of translational orbits of Delone simplices in Delone triangulations of
$d$-dimensional lattices. Let us denote this number $\dt(d)$. 
Since the sum of volumes of representatives from each orbit must equal $d!$, the upper bound 
for $\mv(d)$ implies 
\[
\dt(d) \geq \binom{2d}{d} 2^{-d} = \Omega(2^d / \sqrt{d}).
\]
For an upper bound, we can use essentially the same 
trick as in the proof of Theorem 1  to obtain:

\begin{proposition}
There is a constant $c<1$ such that $\dt(d)\le c^d d!$, for every $d\ge 5$.
\end{proposition}

\begin{proof}
Start with Delone triangulations with the minimal number of 
translational orbits $\dt(d_1)$ and $\dt(d_2)$
in dimensions $d_1$ and $d_2$, and consider the product. 
As before, a perturbation of the product lattice produces a Delone triangulation 
with $\binom{d_1 + d_2}{d_1} \dt(d_1) \dt(d_2)$ translational orbits, from which
we deduce
\[
\frac{\dt(d_1+d_2)}{(d_1+d_2)!} \leq 
\frac{\dt(d_1)}{d_1!} \cdot
\frac{\dt(d_2)}{d_2!}.
\]
Hence,
\[
\left(\frac{\dt(d_1+d_2)}{(d_1+d_2)!}\right)^{\frac{1}{d_1+d_2}} \leq 
\max\left\{\left(\frac{\dt(d_1)}{d_1!}\right)^{\frac{1}{d_1}} ,
\left(\frac{\dt(d_2)}{d_2!}\right)^{\frac{1}{d_2}}\right \}.
\]

Since in every dimension $d\ge 5$ there are Delone simplices of volume greater than one,
we have constants $c_d<1$ such that $\dt(d) = {c_d}^d d!$. It suffices now to take 
as the global constant $c$ for the statement the minimum of $c_5$, $c_6$, $c_7$, $c_8$ and $c_9$.
\end{proof}

Of course, in order to get good (asymptotic) values of the constant in the statement,
one needs to compute the number of translational orbits of Delone simplices in ``good
lattices''.

\item A similar problem that has attracted attention
both for theoretical and applied reasons is the determination of the minimum
number of simplices in a triangulation of the $d$-dimensional cube $[0,1]^d$. 
In fact, our proof of Theorem 1, as well as the upper bound in the previous remark,
is essentially an adaptation 
of the technique used by Haiman \cite{haiman-1991} to construct
``simple and relatively efficient triangulations of the $n$-cube''.

The situation for this problem is that 
the best lower bound known for the minimum size $\operatorname{t}(d)$
of a triangulation of $[0,1]^d$ is
\[
\sqrt{\frac{6}{d+1}} \sim \frac{\sqrt{6}}{\sqrt[d]{2}(d+1)^{\frac{d+1}{2d}} }
\leq \sqrt[d]{\frac{\operatorname{t}(d)}{d!}},
\]
obtained by Smith \cite{smith-2000}. The best (asymptotic) upper bound is
\[
\lim_{d\to \infty} \sqrt[d]{\frac{\operatorname{t}(d)}{d!}} \leq 0.816,
\]
due to Orden and Santos \cite{os-2003}.
For more information on this topic see~\cite{zong-2005}.
\end{enumerate}

\section*{Acknowledgement}

We would like to thank the anonymous referee and Mathieu Dutour Sikiric 
for valuable comments.

\providecommand{\bysame}{\leavevmode\hbox to3em{\hrulefill}\thinspace}
\providecommand{\href}[2]{#2}


\begin{thebibliography}{Dut04}

\bibitem[Bar73]{baranovskii-1973}
E.P. Baranovskii, \emph{Volumes of {L}-simplexes of five-dimensional lattices},
  Math. Notes \textbf{13} (1973), 460--466, translation from Mat. Zametki 13,
  771--782 (1973).

\bibitem[BR98]{br-1998} E.P. Baranovskii and S.S. Ryshkov,
\emph{Repartitioning complexes in $n$-dimensional lattices (with full
description for $n \leq 6$)}, pages 115--124 in Voronoi's Impact on
Modern Science, Book II (P. Engel, H. Syta, eds.; Institute of Math.,
Kyiv 1998 = Vol.21 of Proc. Inst. Math. Nat. Acad. Sci. Ukraine).

\bibitem[CS88]{cs-1988}
J.H. Conway and N.J.A. Sloane, \emph{Sphere packings, lattices and groups},
  Springer-Verlag, New York, 1988.

\bibitem[Del37]{delone-1937}
B.N. Delone, \emph{The geometry of positive quadratic forms}, Uspekhi Mat. Nauk
  \textbf{3} (1937), 16--62, in {R}ussian.

\bibitem[DL97]{dl-1997}
M.M. Deza and M.~Laurent, \emph{Geometry of cuts and metrics}, Springer-Verlag,
  Berlin, 1997.

\bibitem[Dut04]{dutour-2004}
M.~Dutour, \emph{The six-dimensional {D}elaunay polytopes}, European J. Combin.
  \textbf{25} (2004), 535--548.

\bibitem[EK05]{ek05} S. Eliahoua,  and M. Kervaire,
  \emph{A survey on modular Hadamard matrices},
  Discrete Mathematics  \textbf{302}:1--3 (2005), 85--106.

\bibitem[ER02]{er-2002}
R.~Erdahl and K.~Rybnikov, \emph{An infinite series of perfect quadratic forms
  and big {D}elaunay simplices in {$\mathbb Z\sp n$}}, Tr. Mat. Inst. Steklova
  \textbf{239} (2002), 170--178.

\bibitem[GL87]{gl-1987}
P.M. Gruber and C.G. Lekkerkerker, \emph{Geometry of numbers}, North--Holland,
  Amsterdam, 1987.

\bibitem[Hai91]{haiman-1991} 
H. Haiman, \emph{A simple and relatively efficient triangulation
of the $n$-cube}, Discrete Comput. Geom. \textbf{6} (1991), 287--289.

\bibitem[OS03]{os-2003} D. Orden and F. Santos, 
\emph{Asymptotically efficient triangulations of the $d$-cube}, 
Discrete Comput. Geom., \textbf{30}:4 (2003), 509--528.

\bibitem[RS57]{rs-1957}
C.A. Rogers and G.C. Shephard, \emph{The difference body of a convex body},
  Arch. Math. \textbf{8} (1957), 220--233.

\bibitem[Rys76]{ryshkov-1973c}
S.S. Ryshkov, \emph{The perfect form ${A}^k_n$: the existence of lattices with
  a nonfundamental division simplex, and the existence of perfect forms which
  are not {M}inkowski-reducible to forms having identical diagonal
  coefficients}, J. Sov. Math. \textbf{6} (1976), 672--676, translation from
  Zap. Nauchn. Semin. Leningr. Otd. Mat. Inst. Steklova 33, 65--71 (1973).

\bibitem[SV05]{sv-2004}
A.~Sch\"urmann and F.~Vallentin, \emph{Computational approaches to lattice
  packing and covering problems}, Discrete Comp. Geom. (2005), to appear. 
  cf. \href{http://www.arxiv.org/math/0403272}{arXiv:math.MG/0403272}.

\bibitem[Smi00]{smith-2000}
W. D. Smith, \emph{A lower bound for the simplexity of the $N$-cube via hyperbolic volumes}, in ÒCombinatorics of convex polytopesÓ (K. Fukuda and G. M. Ziegler, eds.), European J. Combin., \textbf{21}:1 (2000), 131--137. 

\bibitem[Vor08]{voronoi-1908}
G.F. Voronoi, \emph{Nouvelles applications des param\`etres continus \`a la
  th\'eorie des formes quadratiques. {D}euxi\'eme {M}\'emoire. recherches sur
  les parall\'elloedres primitifs.}, J. Reine Angew. Math. \textbf{134} (1908),
  198--287, and \textbf{136} (1909), 67--181.

\bibitem[Zon05]{zong-2005}
C.~Zong, \emph{What is known about unit cubes}, Bull. of the Amer. Math. Soc.
\textbf{42}:2 (2005), 181--211.

\end{thebibliography}
\end{document}